\documentclass[11pt]{amsart}

\usepackage{amssymb}

\usepackage{graphicx}

\usepackage{amsmath}

\usepackage{verbatim}

\makeindex

\input xy
\xyoption{all}

\newtheorem{theorem}{Theorem}[section]

\newtheorem{corollary}[theorem]{Corollary}

\newtheorem{proposition}[theorem]{Proposition}
\newtheorem{remark}[theorem]{Remark}

\newcommand{\prskip}{\vspace{8pt}} 
\newcommand{\thmskip}{\vspace{10pt}} 
\newcommand{\pfskip}{\vspace{6pt}} 
\newcommand{\sectskip}{\vspace{50pt}} 
\newcommand{\introskip}{\vspace{25pt}} 

\DeclareMathOperator{\diag}{diag}
\DeclareMathOperator{\rank}{rank}

\title{WEYR STRUCTURES OF  MATRICES AND RELEVANCE \\
 TO COMMUTATIVE FINITE-DIMENSIONAL ALGEBRAS} 

\author{ K.\ C.\ O'Meara \\
Department of Mathematics \\
University of Canterbury \\
Christchurch, New Zealand \\  
}
\email{staf198@uclive.ac.nz}

\author{J.\ Watanabe \\
Department of Mathematics \\
Tokai University \\
Hiratsuka 259-1292 Japan \\
Phone +81-463-58-1211 \\
Fax   +81-463-58-9543 \\
Email watanabe.junzo@tokai-u.jp 
}
\email{watanabe.junzo@tokai-u.jp}



\date{March 13, 2017}

\begin{document}

\maketitle  


\begin{abstract}
We relate the Weyr structure of a square matrix $B$ to that of the $t \times t$  block upper triangular matrix $C$ that has $B$ down the main diagonal and first superdiagonal, and zeros elsewhere. Of special interest is the case $t = 2$ and where $C$ is the $n$\,th Sierpinski matrix $B_n$, which is defined inductively by $B_0 = 1$ and $B_n = \left[\begin{array}{cc} B_{n-1} & B_{n-1} \\ 0 & B_{n-1} \end{array} \right]$. This yields an easy derivation of the Weyr structure of $B_n$ as the binomial coefficients arranged in decreasing order.  Earlier proofs of the Jordan analogue of this had often relied on deep theorems from such areas as algebraic geometry. The result has interesting consequences for commutative, finite-dimension algebras.
\end{abstract}

\introskip

Noncommutative finite-dimensional algebras over a field $F$ have been studied, almost without pause, since the  1840's, with many beautiful results uncovered. A sizable group continues to work on them. But interest in their commutative cousins has only recently been revived. The latter study has been less concerned with the intricacies arising from the particular field $F$ than with the radical. In fact, we will assume $F$ is algebraically closed. Of course, finite-dimensional algebras are finitely generated as algebras. For commutative, finite-dimensional algebras $R$ there are a number of simply-stated, basic problems that  remain unanswered. To name one, although we will not pursue this, what is the minimum number of generators of $R$ required in order for some faithful $R$-module $M$ to have dimension (over $F$) less than $\dim R$? If $R$ can be generated by $k$ elements, then for $k = 1,2$, the  minimum dimension of a faithful module is $\dim R$.  For $k >3$ there are easy examples where $\dim M < \dim R$.  But when $k = 3$ this has been open for over 50 years. The question is better known in the form of whether Gerstenhaber's theorem for two commuting $n \times n$ matrices over $F$ also holds for three: if $A,B,C$ are commuting $n \times n$ matrices, must the dimension of the (unital) subalgebra $F[A,B,C]$ of $M_n(F)$ generated by $A,B,C$ have dimension at most $n$? (For those interested in further details, such as how algebraic geometry impacts the problem, see \cite{Ge}, \cite{Gu}, Chapters 5, 7 of \cite{ATLA}, and \cite{HO}. This is another instance where the Weyr form seems  better suited  than its Jordan counterpart.)
\prskip

The application of our theorems on the Weyr structures of block matrices is to the monomial complete intersection ring
\[
    B \ = \ F[x_1,x_2,\ldots,x_n]/(x_1^{d_1+1}, \ldots, x_n^{d_n+1})
\]
where $d_1,d_2,\ldots,d_n$ are integers. If $B_i$ is the homogeneous space of $B$ of degree $i$, we give a relatively simple proof of the result that the multiplication map
\[
\times(x_1 + x_2 + \cdots + x_n)^{N-2k}: B_k \rightarrow B_{N-k}
\]
is a bijection, where $N=d_1 + \cdots + d_n$. This had first been proved by R. Stanley  using the Hard Leftschetz theorem in algebraic geometry, and later by the second author using the theory of the Lie algebra $sl\mbox{(2)}$. A corollary is that the Weyr structure of the multiplication map $\times(x_1 + x_2 + \cdots + x_n): B \rightarrow B$ is the partition $\dim B = \dim B_0 + \dim B_1 + \cdots + \dim B_N$, once the terms are arranged in decreasing order. As background (we will not pursue this connection), the strong Lefschetz property (which can be defined for an endomorphism of any finite graded vector space) is important when interpreted in terms of a representation of the Lie algebra $sl\mbox{(2)}$. The foundation for these representations, in turn, relies on the  Clebsch--Gordan decomposition of modules over $sl\mbox{(2)}$. We mention this in passing because it was discovered by a physicist, and used in quantum mechanics, yet again a reminder that mathematics associated with physics  invariably turns out to be important in other areas.
\prskip

The Weyr structure of certain blocked matrices relates nicely to the Hilbert function of an Artinian algebra, much nicer than the previously known connections in term of Jordan structure. For example, the Hilbert function of the algebra $F[x_1,x_2,\ldots,x_n]/(x_1^2,x_2^2,\ldots,x_n^2)$  has coefficients, when arranged in decreasing order, those in the Weyr structure of the linear map induced by multiplication by a ``general element'' $(1+x_1)(1 + x_2) \cdots (1+ x_n)$. The matrix of the map is, in fact, the $n$\,th Sierpinski matrix $B_n$ described in the \emph{Abstract}.
\prskip

In our Preliminaries section, we record the basic facts about the Weyr form. The latter form has re-emerged in recent years from relative obscurity since its discovery by the Czech mathematician Eduard Weyr in the 1880's. The Weyr form has been shown to a better tool than its Jordan cousin in a number of situations, but the two forms should really be regarded as partners. The user should be  prepared to flip back and forth, using a lovely duality, according to varying situations.
\sectskip

%
%

\section{Preliminaries}
\introskip

The crux of our arguments involve the shifting effect under multiplication by a nilpotent Weyr matrix $W$, as well as the fact that every square matrix $A$ over an algebraically closed field $F$ has a Weyr canonical form. That is, $A$ is similar to a unique Weyr matrix $W \,= \, \diag(W(\lambda_1), \,  W(\lambda_2), \, \ldots, \, W(\lambda_k))$, where $\lambda_1, \lambda_2, \ldots,\lambda_k$ are the distinct eigenvalues of $A$, and where for a given $\lambda \in F$, the \textbf{basic $n \times n$ Weyr matrix $W(\lambda)$} with eigenvalue $\lambda$ takes the following form:
\prskip

\noindent \emph{There is a partition $n_1 + n_2 + \cdots + n_r = n$ of $n$ with $n_1 \ge n_2 \ge \cdots
\ge n_r \ge 1$ such that, when $W(\lambda)$ is viewed as an $r \times r$ blocked matrix $(W_{ij})$, where the $(i,j)$ block $W_{ij}$ is an $n_i \times n_j$ matrix, the following three features are present{\rm :}
\begin{enumerate}
  \item The main diagonal blocks $W_{ii}$ are the $n_i \times n_i$ scalar matrices $\lambda I$ for \linebreak $i=1,\ldots,r$.
  \item The first superdiagonal blocks $W_{i,i+1}$ are full column-rank $n_i \times n_{i+1}$ \linebreak matrices in reduced row-echelon form \ {\rm (}that is, an identity matrix followed by zero rows{\rm )} \ for $i=1,\ldots,r-1$.
  \item    All other blocks of $W$ are zero \ {\rm (}that is, $W_{ij} = 0$ when $j \ne i,i+1${\rm )}.
\end{enumerate}
In this case, we say that $W(\lambda)$ has  \textbf{Weyr structure} $(n_1,n_2,\ldots,n_r)$. If $n_1 = n_2 = \cdots = n_r$, then $W(\lambda)$ is said to have a \textbf{homogeneous} structure.}
\prskip

\noindent For instance, the basic Weyr matrix with eigenvalue $\lambda$ and Weyr structure $(3,3,2,2)$ is

\[
W(\lambda) \ = \ \renewcommand{\arraystretch}{1.15}\left[ \begin{array}{ccc | ccc | cc | cc}
   \lambda  & 0         & 0         & 1         & 0       & 0       &     \multicolumn{4}{c}{}    \\
              & \lambda & 0         & 0         & 1       & 0       &     \multicolumn{4}{c}{}    \\
              &           & \lambda & 0         & 0       & 1       &     \multicolumn{4}{c}{}    \\ \cline{1-8}
              &           &           & \lambda & 0       & 0         & 1   & 0 &   &    \\
              &           &           &           & \lambda & 0       & 0   & 1 &   &    \\
              &           &           &           &        & \lambda  & 0   & 0 &   &    \\ \cline{4-10}
    \multicolumn{6}{c|}{ }                                        & \lambda & 0 & 1 & 0  \\
    \multicolumn{6}{c|}{ }                                        &  & \lambda & 0 & 1  \\ \cline{7-10}
    \multicolumn{8}{c|}{ }        & \lambda & 0  \\
    \multicolumn{8}{c|}{ }        &  0      & \lambda
    \end{array}\right].
\]
The monograph \cite{ATLA} gives a comprehensive account of the Weyr form. See also \cite{HJ} and \cite{S}.
\prskip

For a general square matrix $A$, the \textbf{Weyr structure of $A$ associated with an eigenvalue $\lambda$} is the Weyr structure of the basic Weyr block $W(\lambda)$ that occurs in the unique Weyr form of $A$. This can be calculated without constructing the Weyr form by looking at the ranks (equivalently nullities) of the powers of $A$ (hence much easier than obtaining directly the Jordan structure in terms of ranks of powers, although the Jordan structure can be deduced as the dual partition of the Weyr structure; see \cite{ATLA}, Theorem 2.4.1 and Corollary 2.4.6).
\prskip

\begin{proposition} \label{P:struc} The Weyr structure of a matrix $A$ associated with an eigenvalue $\lambda$ is $(n_1,n_2,\ldots,n_r)$ where

\begin{align}
\ \ \ \ \ r \ &= \  \mbox{nilpotent index of}   \ \  A - \lambda I, \notag \\
          n_i \ &= \  \rank (A - \lambda I) ^{i-1} \ - \ \ \rank (A - \lambda I)^i, \notag
\end{align}

for $i = 1,\ldots,r$. Moreover, $\rank A^i = n_{i+1} + n_{i+2} + \cdots + n_r$ for each $i$.
\end{proposition}
\pfskip

\begin{proof} See \cite{ATLA}, Proposition 2.2.3.
\end{proof}
\prskip

A nilpotent Weyr matrix $W$, of Weyr structure $(n_1,n_2,\ldots,n_r)$,  when \textbf{right} multiplying  a matrix $X$ that is blocked according to the block structure of $W$ (so the $(i,j)$ block $X_{ij}$ is $n_i \times n_j$) shifts the blocks of $X$ one step to the right, introducing a zero first column of blocks and killing the last column of blocks. However, if the Weyr structure of $W$ is nonhomogeneous (meaning not all the $n_i$ are equal),  $W$ can't faithfully shift the $j$th column of blocks of $X$ to the $(j+1)$th column if $n_j > n_{j+1}$.  In this case only the first $n_{j+1}$ columns of $X_{ij}$ are shifted, and the remaining $n_j - n_{j+1}$ are deleted. \textbf{Left} multiplication by $W$ has a similar shifting effect on the rows of blocks of $X$, shifting from the bottom upwards, and appending $n_i - n_{i+1}$ zero rows to $X_{(i+1)j}$  whenever $n_i > n_{i+1}$. See Remark 2.3.1 in \cite{ATLA}.
\prskip

It is critical to our later arguments to have a clear mental picture of what the powers $W^k$ of a nilpotent Weyr matrix $W$ look like. This is best done as a (repeated) special case of what the shifting does when we left or right multiply a matrix $X$ by $W$. We illustrate this in the case $W$ is the nilpotent Weyr matrix of structure $(3,2,2)$. The leftmost matrix $X$ in the product $XW$ here centralizes $W$, so the right hand side of the equation also agrees with left multiplying $X$ by $W$ in terms of shifting rows of blocks. (This and other examples we have used have been purloined from \cite{ATLA}, Chapter 2. But one co-author of the monograph, after consultation with one co-author of the present paper, has agreed not to pursue charges!)

\[
\!\renewcommand{\arraystretch}{1.1}
  \left[ \begin{array}{c c c |c c |c c}
    a & b & e & h & i & l & m \\
    c & d & f & j & k & n & p \\
    0 & 0 & g & 0 & 0 & q & r \\  \cline{1-7}
      &   &   & a & b & h & i \\
      &   &   & c & d & j & k \\   \cline{4-7}
    \multicolumn{5}{c |}{ } & a & b  \\
    \multicolumn{5}{c |}{ } & c & d
    \end{array} \right]\!\!\!
\left[ \begin{array}{c c c |c c |c c}
     0 & 0 & 0 & 1 & 0 & 0 & 0 \\
    0 &  0 & 0 & 0 & 1 & 0 & 0 \\
    0 & 0 &  0 & 0 & 0 & 0 & 0 \\  \cline{1-7}
      &   &   & 0 & 0 & 1 & 0 \\
      &   &   & 0 & 0 & 0 & 1 \\   \cline{4-7}
    \multicolumn{5}{c |}{ } & 0 & 0  \\
    \multicolumn{5}{c |}{ } & 0 & 0
    \end{array} \right]
 =
\left[ \begin{array}{c c c |c c |c c}
    0 & 0 & 0 & a & b & h & i \\
    0 & 0 & 0 & c & d & j & k \\
    0 & 0 & 0 & 0 & 0 & 0 & 0 \\  \cline{1-7}
      &   &   & 0 & 0 & a & b \\
      &   &   & 0 & 0 & c & d \\   \cline{4-7}
    \multicolumn{5}{c |}{ } & 0 & 0  \\
    \multicolumn{5}{c |}{ } & 0 & 0
    \end{array} \right].
\]
\prskip

Our applications of Weyr structures for certain blocked matrices are connected with graded algebras. A graded commutative Artinian $F$-algebra $A$, with grading $A =  \bigoplus_{k = 0}^{n}$, is said to have the \textbf{strong Lefschetz property} if there is a linear element $l$ such that the multiplication map
\[
    \times l^i : A_k \longrightarrow A_{k+i}
\]
has full rank for all $i = 0,\ldots,n$ and $k = 1,\ldots, n-i$. And $A$ has the \textbf{weak Lefschetz property} if the above property holds for all $k$ and $i = 1$. In each case, the element $l$ is referred to as a ``Lefschetz element''. As the second author had observed (before he became aware of the Weyr form and had been using the Jordan form), $A$ having the strong Lefschetz property is equivalent to the Jordan form of the multiplication map by a ``general element'' having Jordan structure the dual partition of $\dim A$ as the sequence of the dimensions of the homogeneous spaces $A_i$, arranged in decreasing order. Therefore, by the duality between the Jordan and Weyr forms, the Weyr structure of the multiplication map must be the sequence of these dimensions (again in decreasing order). The second author also observed  that the weak Lefschetz property is equivalent to the number of blocks in the Jordan form (the nullity of the map) being the maximum $\dim A_i$, and therefore the latter is the size of the first Weyr structure component.
\sectskip

%
%

\section{Block upper triangular $t \times t$ matrices: Reduction to a single eigenvalue of $0$ or $1$}
\introskip

Let $B$ be an $n \times n$ matrix over an algebraically closed field $F$ and consider the $t \times t$ block upper triangular matrix
\[
C \ = \ \renewcommand{\arraystretch}{1.5}\left[\begin{array}{cccccc}
B & B & 0 &  \hdots & 0 & 0 \\
0 & B & B &  \hdots & 0 & 0\\
\vdots & & & & &  \\
0 & 0 & 0 & \hdots & B & B \\
0 & 0 & 0 & \hdots & 0 &  B
\end{array}\right].
\]
By the Generalized Eigenspace Decomposition of $B$ we know $B$ is similar to a block diagonal matrix $\mbox{diag}(B_1,B_2, \ldots, B_k)$ where $B_i$ has a single eigenvalue $\lambda_i$ \,(and so $\lambda_1, \lambda_2,\ldots,\lambda_k$ are the  distinct eigenvalues of $B$). (See, for instance, \cite{ATLA}, Theorem 1.5.2 and Corollary 1.5.4.) If conjugating by an invertible matrix $P$ achieves this decomposition, then it is easily seen that conjugating $C$ by $\mbox{diag}(P,P,\ldots,P)$, followed by conjugation by permutation matrices corresponding to  various transpositions, gives $C = \diag(C_1,C_2,\ldots,C_k)$ where
\[
C_i \ = \ \renewcommand{\arraystretch}{1.5}\left[\begin{array}{cccccc}
B_i & B_i & 0 &  \hdots & 0 & 0 \\
0 & B_i & B_i &  \hdots & 0 & 0\\
\vdots & & & & &  \\
0 & 0 & 0 & \hdots & B_i & B_i \\
0 & 0 & 0 & \hdots & 0 &  B_i
\end{array}\right].
\]

\noindent Since the Weyr structure of $C$ associated with the eigenvalue $\lambda_i$ is the Weyr structure of $C_i$ associated with $\lambda_i$ (which in turn is the Weyr structure of the nilpotent matrix $C_i - \lambda_i I$), the upshot of all this is that it is enough to establish the Weyr structure of $B$ in the case $B$ has a single eigenvalue $\lambda$, that is, $B = \lambda I + W$ where $W$ is a nilpotent matrix. Moreover, since every square matrix is similar to a  Weyr matrix, we can assume $W$ is in fact a nilpotent Weyr matrix because if $Q$ is an invertible matrix such that $Q^{-1}BQ$ is in Weyr form, then conjugating $C$ by $\diag(Q,Q,\ldots,Q)$ gives
\[
C \ = \ \renewcommand{\arraystretch}{2.0}\left[\begin{array}{cccccc}
Q^{-1}BQ & Q^{-1}BQ & 0 &  \hdots & 0 & 0 \\
0 & Q^{-1}BQ & Q^{-1}BQ &  \hdots & 0 & 0\\
\vdots & & & & &  \\
0 & 0 & 0 & \hdots & Q^{-1}BQ & Q^{-1}BQ \\
0 & 0 & 0 & \hdots & 0 &  Q^{-1}BQ
\end{array}\right].
\]
and $Q^{-1}BQ = \lambda I + Q^{-1}WQ$ where $Q^{-1}WQ$ is a nilpotent Weyr matrix.
\prskip

\noindent
\textbf{Henceforth, we assume $B = \lambda I + W$ where $W$ is a nilpotent Weyr matrix with Weyr structure $(m_1,m_2, \ldots, m_r)$.} Note that the Weyr structures of $B$ and $C$ are the same as those of the nilpotent $W = B -\lambda I$ and $X = C - \lambda I$, respectively.
\prskip

If $\lambda = 0$, the connection between the Weyr structures of $B$ and $C$ is easy.  We record this in our next proposition.
\prskip

\begin{proposition}\label{P:WeyrBC}
When $B$ is an $n \times n$ nilpotent matrix and $F$ has characteristic 0 or $p > n$, the Weyr structures $(m_1,m_2, \ldots, m_r)$ and $(n_1,n_2,\ldots,n_s)$ of $B$ and $C$ respectively are related by
\begin{enumerate}
\item
$s = r$.
\item
$n_i = tm_i$ for all $i$.
\end{enumerate}
\end{proposition}\label{P:nilpotent}
\pfskip

\begin{proof} Just looking at the powers of $X^k$ we see that $\rank X^k = t\rank W^k$, whence by Proposition \ref{P:struc}, we have $s = r$ and $n_i = tm_i$ for all $i$.
\end{proof}
\prskip

 When $\lambda \neq 0$ these structures are independent of $\lambda$. This is because a nilpotent matrix and a nonzero scalar multiple of it must be similar --- their powers have the same nullity (see Proposition 2.2.8 in \cite{ATLA}). Hence $X$ is similar to $(1/ \lambda)X$. Also there is an invertible matrix $P$ such that $P^{-1}(1/ \lambda)WP = W$. Therefore conjugating $(1/ \lambda)X$ by $\diag(P,P,\ldots,P)$ shows $X$ is similar to the matrix obtained by replacing $\lambda$ by $1$. \textbf{Henceforth, we can assume $\lambda =  1$ if $\lambda \neq 0$.} To find the Weyr structure of $X$, it is enough by Proposition \ref{P:struc} to find the ranks of powers of $X$.
\sectskip

%
%

\section{The case $t = 2$}
\introskip

\begin{theorem} \label{T:block}
Let $B$ be an $n \times n$ matrix over an algebraically closed field $F$ of characteristic 0 or $p > n$. Let $\lambda$ be an eigenvalue of $B$. Let $C$ be the $2 \times 2$ block upper triangular matrix
\[
    C \ = \  \left[\begin{array}{cc} B & B \\ 0 & B \end{array} \right].
\]
Let $(m_1,m_2, \ldots, m_r)$ and $(n_1,n_2,  \ldots, n_s)$ be the Weyr structures associated with $\lambda$ of $B$ and $C$, respectively. The following relationships hold: \vspace{3mm}

\begin{enumerate}
\item If $\lambda = 0$, then $s = r$ and $n_i = 2m_i$ for all $i$.
\item If $\lambda \neq 0$ and $r = 1$, then $s = 2$ and $n_1 = n_2 = m_1$.
\item If $\lambda \neq 0$ and $r > 1$, then $s = r + 1$ and
\begin{enumerate}
\item  $n_1 = m_1 + m_2$,
\item  $n_{s-1} = m_{s-2}$, $n_s = m_{s-1}$,
\item  $n_i  =  m_{i-1}  + m_{i+1}$ \, for \,$2  \le  i  \le  s-2$.
\end{enumerate}
\end{enumerate}
\end{theorem}
\prskip

\begin{remark} We don't need the algebraically closed assumption if we know a particular $\lambda \in F$ is an eigenvalue of $B$, but in general we do in order to get the reduction in Section 2. \hfill $\square$
\end{remark}
\pfskip

\begin{proof} By our earlier reduction, we can assume $B = \lambda I + W$ where $W$ is a nilpotent Weyr matrix and $\lambda = 0$ or $\lambda = 1$.
\prskip

\noindent CASE (1): $\lambda = 0$.  This is covered by Proposition \ref{P:WeyrBC}.
\prskip

\noindent CASE (2): $\lambda = 1$, $r = 1$. Here $B =  I$, where $I$ denotes the identity matrix of the appropriate size (here $m_1 \times m_1$). Thus $C -  I$ has nilpotent index 2, whence $s = 2$. Also $\rank (C -  I) = m_1$ and so $n_1 = \mbox{nullity} (C - I) = 2m_1 - m_1 = m_1$ and $n_2 = \rank (C -  I) - \rank (C -  I)^2 = m_1$.
\prskip

\noindent CASE (3): $\lambda = 1, r > 1$.  The Weyr structures of $B$ and $C$ are the same as the Weyr structures of $W$ and   \,$X = C - \diag( I, I) = \left[\begin{array}{cc} W & B \\ 0 & W \end{array} \right]$. \,For each integer $k \ge 0$, we have that
\begin{equation} \tag{4}
 X^k = \left[\begin{array}{cc} W^k & kW^{k-1}B \\ 0 & W^k \end{array} \right],
\end{equation}
 and hence the nilpotent index of $X$ is 1 more than that of $W$ (because $B$ is invertible). Therefore $s = r + 1.$
\prskip

We next relate the ranks of $X^k$ and $W^k$. Notice from equation (4) that  $X^k$ is row equivalent to (and hence has the same rank as)
 \begin{equation}\tag{5}
Y \ = \ \left[\begin{array}{cc} W^k & kW^{k-1}B - kW^k \\ 0 & W^k \end{array} \right] \ = \
\left[\begin{array}{cc} W^k & k W^{k-1} \\ 0 & W^k \end{array} \right].
\end{equation}

\noindent
CLAIM:  \emph{for $1 \le k \le s-1$, we have
\[
\rank X^k \ =  \ \rank W^{k-1} \, + \, \rank W^{k+1}.
\]
}
\prskip

\noindent When $k = s-1$, we see directly that $\rank X^k = \rank W^{k-1}$ because the diagonal of $X^k$ is zero. Since $\rank W^{k+1} = 0$, the desired relationship holds. Now assume $1 \le k \le s-2$.  It is enough to show $Y$ has the stated rank. However,  for any nilpotent Weyr matrix W, we can see that the rank of  $Y$ is indeed as claimed, simply by looking at the single nonzero superdiagonal of  blocks in the (1,1), (1,2), and (2,2) blocks of $Y$  (latter of same size as the matrix $W$).  By row operations the (2,2) superdiagonal can be used to clear out the (1,2) superdiagonal  except for the first block which has size $m_1 \times m_k$ and has the identity matrix $I_{m_k}$ as its top half and zeros below. Now all the nonzero rows of $Y$ are independent. At first glance it looks like the rank of $Y$ is therefore $2\rank W^k$. But remember what happens in
the powers of $W$ when $W$ has a nonhomogeneous structure --- blocks are pushed to the right but the last few columns
are lost in a block if there is a squeeze (see discussion at the end of Section 1). Hence we have picked up an extra $m_k - m_{k+1}$ nonzero rows over those in $\diag(W^k,W^k)$. So therefore, using Proposition \ref{P:struc}, we have
\begin{align}
  \rank Y \ &=  \ (m_k - m_{k+1})\, + \, 2\rank W^k \notag \\
            &=  \ (m_k - m_{k+1})\, +  \, 2(m_{k+1} + m_{k+2} + \cdots + m_r) \notag \\
            &= \  (m_k + \cdots + m_r) \, + \, (m_{k+2} + \cdots + m_r) \notag \\
            &= \   \rank W^{k-1} \, + \, \rank W^{k+1}. \notag
\end{align}

\noindent [\emph{To make the above argument clearer, here is the matrix picture of $Y$ in the case $W$ has Weyr structure $(m_1,m_2,m_3,m_4) = (3,2,1,1)$, $\lambda = 1$, and $k = 2$. The header on the matrix indicates the width of the various blocks of $W$.}
\prskip

\begin{align}
& {\Small \begin{array}{cccccccccccccccc}
    &  &  & m_1 & & & m_2 &  m_3 & m_4 &  & m_1  & & & m_2  & m_3 & m_4
 \end{array} } \notag \\
 Y \ =  \  & \ \    \left[\begin{array}{c|c}
 \renewcommand{\arraystretch}{1.15}
  \begin{array}{ccc | cc | c | c}
   0 & 0 & 0 & 0 & 0 & 1 & 0    \\
   0 & 0 & 0 & 0 & 0 & 0 & 0    \\
   0 & 0 & 0 & 0 & 0 & 0 & 0    \\ \cline{1-7}
   \multicolumn{3}{c|}{ }& 0 & 0 & 0 & 1    \\
   \multicolumn{3}{c|}{ }& 0 & 0 & 0 & 0    \\ \cline{4-7}
   \multicolumn{5}{c|}{ } & 0 & 0    \\ \cline{6-7}
   \multicolumn{6}{c|}{ } & 0
       \end{array}    &
   \renewcommand{\arraystretch}{1.15}
  \begin{array}{ccc | cc | c | c}
   0 & 0 & 0 & 2 & 0 & 0 & 0    \\
   0 & 0 & 0 & 0 & 2 & 0 & 0    \\
   0 & 0 & 0 & 0 & 0 & 0 & 0    \\ \cline{1-7}
   \multicolumn{3}{c|}{ }& 0 & 0 & 2 & 0    \\
   \multicolumn{3}{c|}{ }& 0 & 0 & 0 & 0    \\ \cline{4-7}
   \multicolumn{5}{c|}{ } & 0 & 2    \\ \cline{6-7}
   \multicolumn{6}{c|}{ } & 0
       \end{array}   \\ \hline
   \renewcommand{\arraystretch}{1.15}
  \begin{array}{ccccccc}
  \multicolumn{7}{c}{ } \\
  \multicolumn{7}{c}{ } \\
  \multicolumn{7}{c}{ } \\
  \multicolumn{7}{c}{ } \\
  \multicolumn{7}{c}{ } \\
  \multicolumn{7}{c}{ } \\
  \multicolumn{7}{c}{ }
  \end{array}   &
  \renewcommand{\arraystretch}{1.15}
  \begin{array}{ccc | cc | c | c}
   0 & 0 & 0 & 0 & 0 & 1 & 0    \\
   0 & 0 & 0 & 0 & 0 & 0 & 0    \\
   0 & 0 & 0 & 0 & 0 & 0 & 0    \\ \cline{1-7}
   \multicolumn{3}{c|}{ }& 0 & 0 & 0 & 1    \\
   \multicolumn{3}{c|}{ }& 0 & 0 & 0 & 0    \\ \cline{4-7}
   \multicolumn{5}{c|}{ } & 0 & 0    \\ \cline{6-7}
   \multicolumn{6}{c|}{ } & 0
       \end{array}
  \end{array} \right]   \notag
\end{align}
\emph{Notice that here $\rank W^2 = 2$,  and for  $D = \diag(W^2,W^2)$,  we have $\rank D = 2(\rank X^2) = 4$.  However $\rank Y = 5$.  After using the $(2,2)$ block of $Y$ to clear out  the two nonzero entries 2 in the last two columns of the $(1,2)$ block of $Y$, the remaining rows of $Y$ are independent. But we have picked up an extra  \, $ m_2 - m_3 = 1$   independent row in $Y$ in additional to the four in $D$ because of the  $(2,12)$ entry 2 in $Y$.}]
\prskip

Now for the proof of 3(a). From the above Claim and the fact that  $X$ is twice as large as $W$, we have
\begin{align}
n_1 \ &= \ \rank X^0 - \rank X \notag \\
&=\  2\rank W^0 - (\rank W^0 + \rank W^2) \notag \\
&= \ \rank W^0 - \rank W^2 \notag \\
&= \ m_1 \, + \, m_2. \notag
\end{align}
\prskip

For the proof of 3(b), we have
\begin{align}
n_{s-1} \ &= \ \rank X^{s-2} \, -  \, \rank X^{s-1} \notag \\
 &=  \ (\rank W^{s-3} + \rank W^{s-1}) - (\rank W^{s-2} + \rank W^s) \notag \\
 &= \ m_{s-2}, \notag \\
n_s \ &=  \  \rank X^{s-1} \, -  \, \rank X^s \ = \rank W^{s-2} + \rank W^s - 0 \ = \ m_{s-1}. \notag
\end{align}

Finally, for the proof of 3(c), by our Claim, for   $2 \le  i  \le  s-2$ we have that
\begin{align}
n_i \ &= \ \rank X^{i-1} \, - \, \rank X^i  \notag \\
      &=  \ (\rank W^{i-2} + \rank W^i) \, - \, (\rank W^{i-1} + W^{i+1}) \notag \\
      &=  \  m_{i-1} + m_{i+1}, \notag
\end{align}
as  desired.
\end{proof}
\prskip

The sequence of ``Sierpinski'' matrices $B_n$ (an informal term, chosen because the $B_n$ look like a Sierpinski triangle (or gasket), often shown as fractal figure) are defined inductively by $B_0 = 1$,
\[
 B_{n+1} = \left[\begin{array}{cc} B_n & B_n \\ 0 & B_n \end{array} \right]
 \]
 for $n = 0,1,2,\ldots \ .$ Thus $B_n$ is $2^n \times 2^n$. Repeated applications of Theorem \ref{T:block} yield:
\pfskip

\begin{corollary} \label{C:Zilpinski}
Over a field $F$ of characteristic 0 or $p > n$, the Weyr structure of the $n$th Sierpinski matrix $B_n$ is the sequence of binomial coefficients $n!/k!(n-k)!$  for $k = 1,\ldots,n$ arranged in decreasing order. {\rm(}Decreasing for us means non-increasing.{\rm)}
\end{corollary}
\pfskip

\begin{proof} The connection between the $n$th and $(n+1$)th sequences of binomial coefficients (arranged in decreasing  order) is exactly as we have in (2) and (3) of Theorem \ref{T:block} for the Weyr structures of $B_n$ and $B_{n+1}$ (and with $\lambda = 1$). They also have the starting point at $n = 0$. Hence the Weyr structures and sequences of binomial coefficients must be the same.
\end{proof}

Thus, the Weyr structures of the first seven Sierpinski matrices are:
\begin{center}
1 \\
1 \ 1 \\
2 \ 1 \ 1 \\
3 \ 3 \  1 \  1 \\
6 \ 4 \ 4  \ 1 \ 1 \\
10 \ 10 \ 5 \ 5 \ 1 \ 1\\
20 \ 15 \ 15 \ 6 \ 6 \ 1 \ 1
\end{center}

\noindent By contrast, taking the dual structures we obtain the Jordan structures of the first seven Sierpinski matrices:
\begin{center}
1 \\
2 \\
3 \ 1 \\
4 \ 2 \ 2 \\
5 \ 3 \ 3 \ \ 3 \ 1 \ 1 \\
6 \ 4 \ 4 \ 4 \ 4 \ 2 \ 2 \ 2 \ 2 \ 2 \\
7 \  5 \ 5 \ 5 \ 5 \ 5 \ 3 \ 3 \ 3 \ 3 \ 3 \ 3 \ 3 \ 3 \ 3 \ 1 \ 1 \ 1 \ 1 \ 1
\end{center}
\prskip

\noindent So no obvious natural pattern connecting the Jordan structures, although it would be possible to write down a messy relationship by using the Weyr pattern and translating via dual partitions. For instance, the first Jordan structure component of $B_n$  will be $n+1$, the length of the Weyr structure, because this is the nilpotent index of $B_n - I$. And the length of the Jordan structure of $B_n$ will be the first Weyr component, because this is the nullity of $B_n - I$. This is why we didn't display the Jordan structure of $B_7$ because it has 35 Jordan structure components (as against 8 Weyr components)! Is this yet another situation where the Weyr form seems more in tune to natural phenomenon than its Jordan counterpart?
\prskip

On the other hand, using the Lie algebra $sl\mbox{(2)}$, the second author \cite{WAT} established a connection between the Jordan structure of the multiplication map by ``a general element'' of  $B = F[x_1,x_2,\ldots,x_n]/(x_1^{e}, x_2^{e}, \ldots, x_n^{e})$ and the $n$\,th Sierpinski matrix $B_n$. He did this by showing that the Jordan structure of the matrix of the  multiplication map of a general element is the dual of the sequence (in decreasing order) of the coefficients  of the Hilbert function, which here is given by
\[
(1+T+ T^2 + \cdots + T^{e-1})^ n.
\]
Using simpler methods, Hidemi Ikeda much later proved the same thing for $B = F[x_1,x_2, \ldots, x_n]  /  (x_1^2, x_2^2, \ldots , x_n^2)$.
\sectskip

%
%

\section{The case $t = 3$}
\introskip

Examined closely, the proof for $t = 2$ is actually very simple. However, it gives little indication of what happens when $t > 2$. Moreover, a critical point in the argument later for $t > 2$ doesn't occur when $t = 2$. The case $t = 3$ is a better indicator of what  happens in general, and the pitfalls to watch out for, but even here one is left guessing the general pattern. It is the case $t = 4$ (combined with $t = 2$ and $t = 3$) that finally strongly suggests the general pattern, as well as the inductive argument to use.
\prskip

\begin{proposition}\label{P:ranks} Again suppose $B$ is an $n \times n$ matrix over a field $F$ of characteristic 0 or $p > n$, and assume $B$ has a single eigenvalue $\lambda$, and this is nonzero. Let $t = 3$ and let $C$ be the $t \times t$ block upper triangular matrix defined earlier. Let $W = B - \lambda I$ and $X = C - \lambda I$. Let $r,s$ be the nilpotent indices of $W, C$ respectively.   Then we have:
\begin{enumerate}
\item
 $s = r + 2$.
\item
$\rank X \, = \, 2n + \,  \rank W^3$.
\item
For $2 \le k \le s - 1$,
\[
\rank X^k \, = \, \rank W^{k-2} \, + \, \rank W^k \, + \, \rank W^{k+2}.
\]
\end{enumerate}
\end{proposition}
\prskip

\begin{proof} Let $(m_1,m_2,\ldots,m_r)$ and $(n_1,n_2,\ldots, n_s)$ be the Weyr structures of $W$ and $X$ respectively. By our earlier reduction in Section 2, we can assume $\lambda = 1$.
\prskip

\noindent (1) and (3). Assume $2 \le k \le r+1$. We have

\[
X^k \, = \, \renewcommand{\arraystretch}{3.0}\left[\begin{array}{ccc}
W^k & \ \ kW^{k-1}B & k(k-1)/2W^{k-2}B^2 \\
0   & \ W^k & kW^{k-1}B \\
0 & 0 & W^k
\end{array}\right].
\]
Hence the nilpotent index of $X$ is 2 more than that of $W$. Thus $s = r+2.$  Expanding the terms by replacing $B$ with $ I + W$, we see that $X^k$ has the form
\[
   \ \renewcommand{\arraystretch}{3.0}\left[\begin{array}{cccc}
     W^k & \ \ aW^{k-1}+ bW^k \ \ & cW^{k-2}+ dW^{k-1} + eW^k  \\
     0 & W^k & aW^{k-1}+ bW^k   \\
     0 & 0 & W^k
     \end{array}\right]
\]
for nonzero integers  $a,b,c,d,e$ and with $a = k$ and  $c = k(k-1)/2$. Using row and column operations, we see that $X^k$ is equivalent to
\[
  Y \ = \ \renewcommand{\arraystretch}{3.0}\left[\begin{array}{cccc}
     W^k & \ \ aW^{k-1} \ \ & bW^{k-2}  \\
     0 & W^k & aW^{k-1}   \\
     0 & 0 & W^k
     \end{array}\right]
\]
for $a = k,\, b = k(k-1)/2$ \, (so $b$ has been renamed as the above $c$).
\prskip

When $k = r+1$, only the $(1,3)$ entry of $Y$ is nonzero, so clearly $\rank Y = \rank W^{k-2}$ and the relation in (3) holds (because $W^k = W^{k-1} = 0)$. Henceforth we assume $2 \le k \le r$.
\prskip

When we refer to the ``blocks of $Y$'' (or of $X^k$) we mean the $9$ blocks resulting from partitioning $3n$ as $(n,n,n)$. But we also need to refer to entries within these blocks, and there it is convenient to reference ``blocks within a block'' by partitioning an $n \times n$ matrix by the Weyr structure $(m_1, m_2, \ldots,m_r)$ of $W$.
\prskip

We will find $\rank Y$ using only column operations, moving across the three columns of blocks, and ensuring the nonzero individual  columns (among the $3n$ columns of $Y$) are independent at each step. This will establish a pattern which will help us with larger  $t$  by induction. The first column of blocks causes no problem because the nonzero columns of $W^k$ are linearly independent and there are $\rank W^k = m_{k+1} + \cdots + m_r$ of them.
\prskip

Now move to the 2nd column of blocks. The columns of $Y$ that contain a  nonzero column of $W^k$ in the $(2,2)$ block are independent, and independent of all columns to their left in $Y$.  So they contribute another $\rank W^k$ to the $\rank Y$. \emph{Therefore, for rank purposes, the only columns in $Y$ that can further contribute must not ``step over the $W^k$ line''.} (Note that  the nonzero columns of $W^k$ begin at the start of the $(1,k+1)$ block, when we block according to the partition $(m_1,m_2,\ldots,m_r)$, and are then ALL nonzero from that point on.) Now the only stuff in the $(1,2)$ block of $Y$ that is inside the $W^k$ line is the $(1,k)$ block of $aW^{k-1}$, which has $aI_{m_k}$ as its top part and zeros below. Clear out the first $m_{k+1}$ columns of $aW^{k-1}$ using those from the $(1,k+1)$ block of $W^k$ in the $(1,1)$ block of $Y$. Now all the nonzero columns within the first two blocks of $Y$ are independent. Thus the column rank of $Y$ to this point is
\[
2\rank W^k \, + \,  (m_k - m_{k+1}).
\]
Moreover, the column space of $Y$ to this point includes the space of all $3n \times 3n$ column matrices with arbitrary entries in the first $m_k$ positions and zeros after (the natural copy of $F^{m_k}$, the space of $m_k \times 1$ column vectors over $F$). This second step has been quite straightforward. Moreover, \emph{we never have to change this part of $Y$ in the third step!}.
\prskip

Now to the third column of blocks.  Again there is a contribution of $\rank W^k$ to $\rank Y$ coming from the $W^k$ in the $(3,3)$ block, and any further contribution can only be within the $W^k$ line. The only stuff within the line, in terms of the blocking of $W$,  are two blocks in the $(1,3)$ block of $X^k$, and one in the $(2,3)$ block. The picture is this:
{\Small
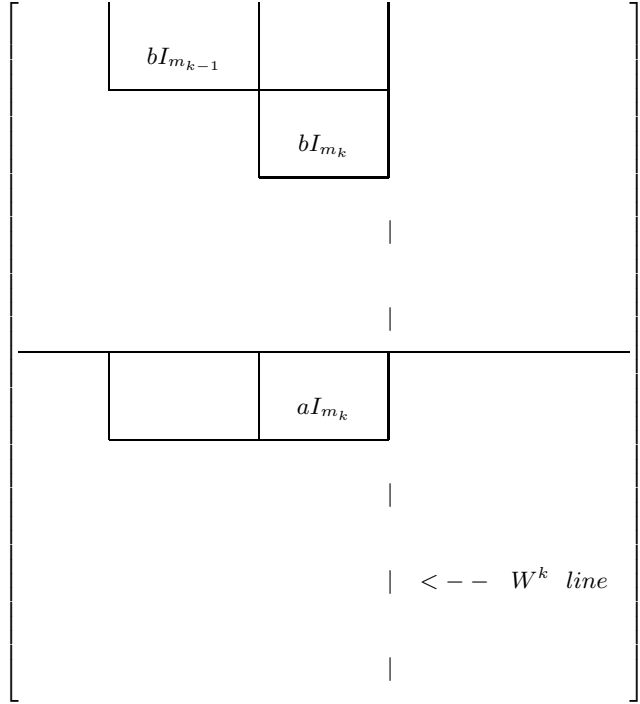
\[
\renewcommand{\arraystretch}{3.0}\left[\begin{array}{cccc |c|c| ccccccc}
\multicolumn{4}{c|}{\phantom{I_{m_{k-1}}}} & \ \ \  bI_{m_{k-1}} \ \ \  &  & \multicolumn{7}{|c}{\phantom{I_{m_{k-1}}, I_{m_{k-1}}, I_{m_{k-1}}}}\\ \cline{5-6}
\multicolumn{5}{c|}{ } & bI_{m_k} & \multicolumn{7}{|c}{ } \\ \cline{6-6}
\multicolumn{5}{c}{ } & \multicolumn{8}{c}{ }
{$\begin{array}{ccccccc}
   | & \multicolumn{6}{c}{ } \\
  | & \multicolumn{6}{c}{\phantom{<-- \ \ W^k \ \ line}} \\ 
    \end{array}$} \\ \hline
\multicolumn{4}{c|}{\phantom{ I_{m_{k-1}}}  } & {\phantom{ \ \ \  bI_{m_{k-1}}} \ \ \   } &  \ \ \ aI_{m_k} \ \ \  & \multicolumn{7}{|c}{\phantom{I_{m_{k-1}}, I_{m_{k-1}}, I_{m_{k-1}}}}\\ \cline{5-6}
\multicolumn{5}{c}{ } & \multicolumn{8}{c}{ }
{$\begin{array}{ccccccc}
   | & \multicolumn{6}{c}{ } \\
  | & \multicolumn{6}{c}{<-- \ \ W^k \ \ line} \\ 
  | & \multicolumn{6}{c}{ } \\
  \end{array}$}
\end{array}\right]
\]
}
Knowing we have the copy of $F^{m_k}$ in the top part of our column space to date, we can clear out the first $m_k$ columns of the $bI_{m_{k-1}}$ in the first of these. What remains from this block contributes another $m_{k-1} - m_k$ to $\rank Y$. This now leaves only the other two blocks containing $bI_{m_k}$ and $aI_{m_k}$, in the $(2,k)$ block and $(1,k)$ block respectively (relative to the structure of $W$, but in the $(1,3)$ and $(2,3)$ blocks of $Y$). \textbf{\emph{But here we must be very careful about the clearing argument involved. It depends on  $b \neq a^2$}} \textbf{!!} Thus any ``hand-waving'' argument based on the matrix picture, but not taking into account the actual values of $a$ and $b$, will likely give the wrong answer! (The first author has got Matlab to thank for rescuing him after falling into this pit by assuming he could work with just $a = b = 1$\,!)
\prskip

To clear out a particular column of blocks in the 3rd column of $Y$, we look to the left and note all columns of blocks that have a common part with the column under consideration. Then we argue what linear combinations of the stuff to the left can be used for clearing. In turn that involves looking at the rank of a small matrix and deciding what is the dependence of its last column on the earlier ones. But this requires a very clear mental picture of what a power of the Weyr matrix $W$ looks like. If $b = a^2$, then subtracting $a$ times the matching blocks in the second column of blocks of $Y$ from the two blocks in question in the third column of blocks gives the desired clearing for new independent columns. This results in an extra contribution of $m_k - m_{k+1}$ to $\rank Y$. But when $b \neq a^2$  the matching blocks in the second column are independent of those in the third,  because the matrix $\left[\begin{array}{cc} a & b \\ 1 & a \end{array}\right]$ is nonsingular. However, we can subtract from those in the third, $a$ times those in the second column, followed by subtracting $b$ times the $k+2$ column of $W^k$ within the $(1,1)$ block of $Y$. In this case, the extra contribution to $\rank Y$ is $m_k - m_{k+2}$. But when we look at the specific $a = k$ and $b = k(k-1)/2$ we are dealing with, we see $b \neq a^2$. Thus it is the latter contribution that applies. Hence we have established that
\begin{align}
\rank X^k \ &= \ 3\rank W^k \, + \, (m_k - m_{k+1}) \, + \, (m_{k-1} - m_k) \, + \, (m_k - m_{k+2}) \notag \\
            &= \ 3\rank W^k \, + m_{k-1} \, + \, m_k \, - \, m_{k+1} \, - \, m_{k+2} \notag \\
            &= \ 3(m_{k+1} + \cdots + m_r) \, + \,   m_{k-1} \,  + \, m_k\,  - \, m_{k+1} \, - \, m_{k+2} \notag \\
            &= \ (m_{k-1} + \cdots + m_r) \, + \, (m_{k+1} + \cdots + m_r) + (m_{k+3} + \cdots + m_r) \notag \\
            &= \ \rank W^{k-2} \, + \, \rank W^k \, + \, \rank W^{k+2}. \notag
\end{align}
\prskip

\noindent (2) When $k = 1$, the $(1,3)$ block of $X^k$ is zero. So the argument is slightly different but simpler.
\end{proof}
\thmskip

\begin{corollary} Let $(m_1,m_2,\ldots,m_r)$ and $(n_1,n_2,\ldots, n_s)$ be the Weyr structures relative to a nonzero eigenvalue $\lambda$ of $B$ and $C$ respectively. We have:
\begin{enumerate}
\item
$s = r + 2$.
\item
$n_1 =  m_1 + m_2 + m_3$.
\item
$n_2 = m_1 + m_2 + m_4$.
\item
$n_i \ = \ m_{i-2} \, + \, m_i \, + \, m_{i+2}$ \, for \, $3 \le i \le s - 3$.
\item
$n_{s-2} = m_{s-4} + m_{s-2}$, \, $n_{s-1} = m_{s-3}$, \, and \, $n_s = m_{s-2}$.
\end{enumerate}
\end{corollary}
\pfskip

\begin{proof} (1) This was established in  Proposition \ref{P:ranks}(1).
\prskip

\noindent (2), (3). By Propositions  \ref{P:struc} and \ref{P:ranks}, we have
\begin{align}
n_1 \  &= \  \rank I - \rank X = 3n - (2n + \rank W^3) \notag \\
       &= \ n - (m_4 + m_5 + \cdots + m_r) = (m_1 + \cdots + m_r) - (m_4 + \cdots + m_r) \notag \\
       &= \ m_1 + m_2 + m_3, \ \ \mbox{and} \notag \\
n_2 \  &=  \ \rank X - \rank X^2 = (2n + \rank W^3) - (\rank W^0 + \rank W^2 + \rank W^4) \notag \\
       &= \ (n - \rank W^2) + (\rank W^3 - \rank W^4) \notag \\
       &= \ m_1 + m_2 + m_4. \notag
\end{align}
\prskip

\noindent (4) For $3 \le i \le s-3$, we have by Proposition \ref{P:ranks}
\begin{align}
n_i \ &= \ \rank X^{i-1} \, - \, \rank X^i \notag \\
      &= \ (\rank W^{i-3} \, + \, \rank W^{i-1} \ + \, \rank W^{i+1}) \, - \, (\rank W^{i-2} \, + \, \rank W^i \ + \, \rank W^{i+2}) \notag \\
      &= \ (\rank W^{i-3} \, - \rank W^{i-2}) \, + \, (\rank W^{i-1} \, - \rank W^i) \, + \, (\rank W^{i+1} \, - \rank W^{i+2}) \notag \\
      &= \ m_{i-2} \, + \, m_i \, + \, m_{i+2}. \notag
\end{align}
\prskip

\noindent (5) Same calculation as in (4) after noting $\rank W^j = 0$ for $j \ge s - 2$.
\end{proof}
\sectskip

%
%

\section{The case $t = 4$ and higher}
\introskip

This involves the same strategy, working out the contribution to $\rank X^k$ that comes from the new column of blocks (the 4th here or column $t$ in general). We have
\[
X^k \ = \ \renewcommand{\arraystretch}{2.5}\left[\begin{array}{cccc}
   W^k & a_kW^{k-1}B & b_kW^{k-2}B^2 & c_kW^{k-3}B^3 \\
   0 & W^k & a_kW^{k-1}B & b_kW^{k-2}B^2  \\
   0 & 0 & W^k & a_kW^{k-1}B  \\
   0 & 0 & 0 & W^k
   \end{array}\right]
\]
where \, $a_k = k$, \, $b_1 = 0$ and $b_k = k(k-1)/2$ for $k \ge 2$, \, $c_1 = c_2 = 0$ and $c_k = (k-2)(k-1)k/6$ for $k\ge 3$. Again by our earlier reduction we can assume $B$ has just one eigenvalue $\lambda$, and that $\lambda = 1$. From the expression for $X^k$ we see that the nilpotent index of $X$ is 3 more than that of $W$, whence $s = r+3$. Expanding terms using $B =  I + W$, and clearing using row and column operations shows $X^k$ is equivalent to
\[
Y \ = \ \renewcommand{\arraystretch}{2.5}\left[\begin{array}{cccc}
   W^k & a_kW^{k-1} & b_kW^{k-2} & c_kW^{k-3} \\
   0 & W^k & a_kW^{k-1} & b_kW^{k-2} + d_kW^{k-1} \\
   0 & 0 & W^k & a_kW^{k-1}  \\
   0 & 0 & 0 & W^k
   \end{array}\right]
\]
where $a_k, b_k, c_k$ are the integers as above. The coefficient $d_k$ is a nonzero integer whose value need not concern us. When $t = 3$, we were able to clear out so as to leave just the smallest power of $W$ in each of the nonzero blocks, but it looks like this is not possible when $t = 4$ (getting rid of $d_kW^{k-1}$ without messing other things up). However, this type of clearing was just for convenience, and the arguments work just as well without it (because a given power of $W$ ``covers all the higher powers'' in the sense that its column space contains the column spaces of higher powers).
\prskip

Suppose $k \ge 3$ (the arguments for $k = 1$ and $k = 2$ are slightly different but easier). Again we can ignore the stuff outside the $W^k$ line. Blocking the 4th column of $Y$ (which is a $4n \times n$ matrix) by partitioning its columns by $(m_1,m_2,\ldots,m_r)$ and its rows by 4 lots of this partition, we see that the matrix inside the $W^k$ line has one (nonzero) block in column $k-2$ (an  $m_1 \times m_{k-2}$ matrix), two  blocks  in column $k-1$, and three  blocks in column $k$. The contribution to $\rank X^k$ from column $k-2$ is $m_{k-2} - m_{k-1}$. The contribution from column $k-1$ is $m_{k-1} - m_{k+1}$. For this we observe that the matrix $\left[\begin{array}{ccc} a_k & b_k & c_k \\ 1 & a_k & b_k \end{array}\right]$ has rank 2 with the the first two columns independent (whence the 3rd column is a linear combination of the first two). The contribution from column $k$ is $m_k - m_{k+3}$, and here we use the fact that
 \[
 \left[\begin{array}{cccc} 1 & a_k & b_k & c_k \\ 0 & 1 & a_k & b_k \\ 0 & 0 & 1 & a_k \end{array}\right]
\]
has rank 3 with the first 3 columns independent (so the 4th column is a combination of the first three), as well the previous observation concerning the $2 \times 3$ matrix. Of course, outside the $W^k$ line we have a contribution of $\rank W^k$. Thus the total contribution from all of the 4th column of $Y$ is
\[
   \rank W^k \, + \, m_{k-2} \, + \, m_k \, - \, m_{k+1} \, - \, m_{k+3}.
\]
Adding to this the known contribution of \, $3 \rank W^k \, + \, m_{k-1} \, + \, m_k \, - \, m_{k+1} \, - \, m_{k+2}$ from the first 3 columns of $X^k$ (the case $t = 3$), we have:
\begin{align}
\rank X^k \  &= \  \rank Y \notag \\
  &= \ 4\rank W^k \, + \, m_{k-2} \, + \, m_{k-1} \, + \, 2m_k \, - \, 2m_{k+1} \, - \, m_{k+2} \, - m_{k+3} \notag \\
             &= \ \rank W^{k-3} \, + \, \rank W^{k-1} \, + \, \rank W^{k+1} \, + \, \rank W^{k+3}. \notag
\end{align}
\prskip

If we also compute \, $\rank X = 4\rank W + 3m_1 - m_2 - m_3 -m_4 = 3n + \rank W^4$\, and \, $\rank X^2 = 4\rank W^2 + 2m_1 + 2m_2 - 2m_3 - m_4 - m_5$, we can compute the new Weyr structure components $n_1, n_2, \ldots n_s$ using the connection $n_i = \rank X^{i-1} - \rank X^i$ as:
\prskip

\begin{proposition} For $t = 4$ we have:
\begin{enumerate}
\item
$s = r + 3$.
\item
$n_1 = m_1 + m_2 + m_3 + m_4$.
\item
$n_2 = m_1 + m_2 + m_3 + m_5$.
\item
$n_3 = m_1 + m_2 + m_4 + m_6$.
\item
$n_i \, = \, m_{i-3} \, + \, m_{i-1} \, + \, m_{i+1} \, + \, m_{i+3}$ \ \ for \, $4 \le i \le s - 4$.
\item
$n_{s-3} = m_{s-6} + m_{s-4},\, n_{s-2} = m_{s-5} + m_{s-3}, \, n_{s-1} = m_{s-4}, \, n_s = m_{s-3}.$
\end{enumerate}
\end{proposition}
\prskip

\emph{\textbf{The pattern is now perfectly clear.}} For a general $t$, we have
\[
 \rank X^k \ = \ \rank W^{k-t+1} \, + \, \rank W^{k-t+3} \, + \, \rank W^{k-t+5} \, + \,  \cdots \, + \rank W^{k+t-1}
\]
for $k \ge t - 1$. Also, $s = r + t - 1$ and the ``middle range'' of  Weyr structure components $n_i$ are given by
\[
   n_i \ = \ m_{i - t + 1} \, + \, m_{i-t + 3} \, + \, m_{i - t + 5} \, + \, \cdots \, + \, m_{i + t -1}
\]
for $t \le i \le s - t$. So the rule is \emph{``go back $t-1$ terms from $m_i$  to get the first term on the right hand side, then include  all terms got by going up in steps of 2, until you have a total of $t$ terms''.} Actually this also holds for all $i \ge t$ if we ignore terms that no longer make sense (such as $m_{r+1}$). The initial $n_i$ for $i = 1,\ldots, t-1$ are given by
\[
n_i \ = \ m_1 \, + \cdots \, + m_{t-i+1} \, + \, m_{t-i+3} \, + \, m_{t-i+5} \, + \, \cdots \, + \, m_{t+i-1}.
\]
So the rule here is to \emph{``add the first $t-i+1$ terms before going up in steps of 2 to reach a total of $t$ terms''.}
\prskip

However, a proof of these claims for general $t$ requires a careful argument. We won't give the full details, just a sketch. For the induction to work in going from $t-1$ to $t$, it is enough to show that the total contribution to $\rank X^k$ from the new $t$\,th column is the proposed $\rank X^k$ for $t$ less the assumed one for $t-1$. As before, there is a contribution of $\rank W^k$ to the right of the $W^k$ line, so we need to know the extra contribution $E$ that comes from left of the $W^k$ line. A little arithmetic shows this must be (for the induction to work)
\[
E \ = \ \left\{\begin{array}{lcl}
\ \sum_{i=1}^{(t-1)/2}\, m_{k-t+2i} & - & \sum_{i=1}^{(t-1)/2}\, m_{k+2i} \ \ \ \ \mbox{if $t$ is odd} \\
                                 \\
     \ \sum_{i=1}^{t/2}\ \ m_{k-t+2i} & - & \sum_{i=1}^{t/2}\ \ m_{k+2i} \ \ \ \ \ \ \ \mbox{if $t$ is even.}
     \end{array}\right.
\]
(Note these expressions contain no redundant $m_j$.) Next, instead of the labels $a_k,b_k,c_k$ we used for the coefficients of the smallest power of $W$ in the top row of blocks of the $4 \times 4$ matrix $X^k$, for general $t$ and fixed $k$, we denote the coefficient of the smallest power of $W$ in the $(1,i)$ block by $a_i$. Thus $a_1 = 1, \, a_2 = k(k-~1)/2, \, a_3 = k(k-1)(k-2)/6, \ldots$ and the general $a_i$ for $i \ge 1$ can be computed inductively to be $a_i = {k \choose i}$. It is important in the arguments that follow to note that, due to the actual values of the $a_i$, each of  the $b \times b$ submatrices of the $b \times t$ matrix
\[
B \ = \ \left[\begin{array}{cccccc}
a_1 & a_2 & a_3 & a_4 & \hdots & a_t \\
0 & a_1 & a_2 & a_3 & \hdots & a_{t-1} \\
\vdots & & & & &  \\
0 & \hdots & a_1 & a_2 & \hdots & a_{t-b+1}
\end{array}\right]
\]
is nonsingular for $b = 2,\ldots t$. In particular, the last column of $B$ is a linear combination of the previous $b$ columns, but no fewer.
\prskip

Blocking the $t$\,th column of $X^k$ (which is a $tn \times t$ matrix) by partitioning its columns by $(m_1,m_2,\dots, m_r)$ and its rows by $t$ lots of this partition, we then calculate what extra contribution each nonzero column of blocks inside the $W^k$ line makes to the already known contribution from columns $1,2,\ldots,t-1$ of the $t \times t$ blocked $X^k$. Using the same sort of the argument we used for $t = 4$ yields the following results (where $\#$ blocks is the number of nonzero  blocks in a particular column of blocks):
\begin{table}[h]
\begin{center}
\renewcommand{\arraystretch}{1.3}
\begin{tabular}{|c|c|c|}
\hline
column & $\#$ blocks & contribution \\ \hline
$k-t+2$ & 1 & $m_{k-t+2} \, - \, m_{k-t+3}$ \\
$k-t+3$ & 2 & $m_{k-t+3} \, - \, m_{k-t+5}$ \\
$k-t+4$ & 3 & $m_{k-t+4} \, - \, m_{k-t+7}$ \\
$k-t+5$ & 4 & $m_{k-t+5} \, - \, m_{k-t+9}$ \\
\vdots & \vdots & \vdots\\
$k$ & $t-1$ & $m_k \, - \, m_{k+t-1}$ \\ \hline
\end{tabular}
\end{center}
\end{table}

\noindent Adding these contributions leads to the above value of $E$.
\sectskip

%
%

\section{Relevance to Commutative Finite-dimensional algebras}
\introskip

Let $F$ be a field and let $R=F[x_1, x_2, \cdots, x_n]$ be the polynomial ring  in $n$ variables with coefficients in $F$. Let $A=R/I$, where $I$ is an ideal which contains
\[
(x_1^{d_1+1}, \, x_2^{d_2+1},  \ldots , \,x_n^{d_n+1})
\]
for some integers $d_1,d_2,\ldots,d_n$. In this case $A$ is a commutative finite-dimensional algebra over $F$.  Any element $f \in A$ induces via multiplication an endomorphism of the vector space $\times f : A \to A$, which has a single eigenvalue (namely, the constant part of $f$) and, in most cases, a complicated nilpotent part. One of the  basic problems in the theory of  Artinian rings is  to determine the Weyr form of $\times f$  for a general element $f \in A$.
\prskip

If $I=(x_1^{d_1+1}, x_2^{d_2+1},  \ldots , x_n^{d_n+1})$, then $A=R/I$ is called a \textbf{monomial complete intersection ring}. First we treat the quadratic monomial complete intersection:
\[
d_1=d_2=\cdots =d_n=1.
\]
For the rest of this section we fix
\[
A=F[x_1, x_2, \ldots, x_n]/(x_1^2, x_2^2, \cdots, x_n^2).
\]
The set of square-free monomials in $x_1, \ldots, x_n$ is a basis for $A$. We order these using the reverse lexicographic order:
\begin{center}
$1 \, < \, x_1 \, < \, x_2 \, < \, x_1x_2 \, < \, x_3 \, < \, x_1x_3 \, < \, x_2x_3 \, < \,  x_1x_2x_3 \, < \, x_4 \, < $ \\
$ \cdots \, < \, x_2x_3 \, \cdots x_n \, < \, x_1x_2 \cdots x_n.$
\end{center}
This sequence can be characterized inductively by saying that the monomials in the second half are divisible by $x_n$, and if we substitute 1 for $x_n$, the sequence coincides with the first half in the same order. We fix as our ordered basis for $A$ the set of square-free monomials ordered in this way. Note that the multiplication by any variable $x_i$ kills a monomial which is divisible by $x_i$ but it preserves the order of the rest of the basis.
\prskip

Let $g$ be the sum of all the square-free monomials, that is,
\begin{align}
g \ &= \ 1+x_1 + x_2 + \cdots + x_1x_2 + x_1x_3 + \cdots + x_1x_2\cdots x_n \notag \\
    &= \ (1+x_1)(1+x_2) \cdots (1+x_n). \notag
\end{align}
We consider the linear map $\times g: A \to A$ defined by $m \mapsto gm$. The matrix for $\times g$ relative to our ordered basis  above is exactly the $n$th Sierpinski $B_n$ described in Section 3. Hence from Corollary \ref{C:Zilpinski} we obtain the following description of the Weyr structure of $\times g$:
\thmskip

\begin{theorem}
Let $p$ be the characteristic of $F$. Assume that either $p=0$  or   $p > n$. Then the  Weyr structure of $\times g$ (associated with its sole eigenvalue 1) is the sequence of  binomial coefficients ${n \choose i}$ arranged in the decreasing order.
\end{theorem}
\prskip

Let $l = x_1+x_2 + \cdots + x_n$, and consider the multiplication map $\times (1+l): A \to A$.  Exactly the same argument as for $\times g$ yields the following:
\thmskip

\begin{theorem}\label{slp_quad_mono_ci}
\begin{enumerate}
\item The Weyr structure of $\times (1+l):A \to A$ is the sequence of binomial coefficients arranged in  decreasing order.
\item The Weyr structure of $\times l:A \to A$ is the sequences of binomial coefficients arranged in  decreasing order.
\end{enumerate}
\end{theorem}
\prskip

Note that the algebra  $A$ has a natural grading and $A$ may be regarded as a graded algebra. In fact, if we denote by $A_k$ the vector space spanned by the square-free monomials of degree $k$,
\[
A_k = \{ x_{j_1}x_{j_2}\cdots x_{j_k} \, : \,1 \leq j_1 < j_2 < \cdots < j_k \leq n  \},
\]
then  $A$  decomposes as a direct sum of subspaces
\[
A=\bigoplus _{k=0} ^n A_k.
\]
Note that the multiplication of A, that is $(f,g) \mapsto fg$, is compatible with the grading, meaning that the multiplication restricts to a map $A_i \times A_j \to A_{i+j}$ for each pair $i,j$. Notice also that $\dim _F A_k = \dim _F A_{n-k} = {n \choose k}$, and the multiplication map $\times l : A \to A$ acts on the summands $A_k$ of our grading by mapping
\[
A_0 \, \stackrel{\times l}{\to} \, A_1 \, \stackrel{\times l}{\to} \, A_2 \, \stackrel{\times l}{\to} \, A_3 \, \, \to \cdots \, \to \, A_{n}.
\]
It is easy to see the following is true:
\prskip

\begin{proposition}
\[
\mbox{\rm rank }\left[\times l^k :A \to A \right] \ = \ \sum _{i=0}^{n-k}\mbox{\rm rank }\left[\times l ^k: A_i \to A_{i+k}\right].
\]
\end{proposition}
\pfskip

By Theorem~\ref{slp_quad_mono_ci} and Proposition \ref{P:struc}, this is possible only if $\times l ^k: A_i \to A_{i+k}$ has  full rank.
We single this out as a theorem:

\begin{theorem} For all $k=0,1,2, \cdots, [n/2]$, the multiplication map
  $\times l^{n-2k}: A_k \to A_{n-k}$ is bijective.
\end{theorem}
\pfskip

Let $n=n_1+n_2+ \cdots + n_r$ be a partition of $n$ (so the  $n_i$ are integers with $n_1 \ge n_2 \ge \cdots \ge n_r > 0$). Divide the  set of variables $x_1, \cdots , x_n$ into $r$ groups
\[
\{ \underbrace{x_1,\ldots, x_{n_1}}_{n_1} \} \ \ \  \{ \underbrace{x_{n_1+1},\ldots, x_{n_1+ n_2}}_{n_2} \} \ \ \ \cdots \ \ \ \{ \underbrace{x_{n_1+ \cdots + n_{r-1}+1}, \ldots, x_n}_{n_r} \}.
\]
Let  $Q=F[y_1, y_2, \cdots, y_r]$ be the polynomial ring in $r$ variables, and  define the ring homomorphism $ \phi :Q \to A$ by sending $y_i$ to the sum of the variables in the $i$-th group:
\[
\phi(y_i)=x_{n_1+n_2 +n_{i-1}+1 } +x_{n_1+n_2 +n_{i-1}+2 } + \cdots + x_{n_1+n_2 + \cdots +n_i}.
\]
It is not difficult to see that the kernel of $\phi$ is the ideal generated by
\[
(y_1^{n_1+1},y_2^{n_2+1}, \ldots, y_r^{n_r+1}).
\]
Hence we have a natural inclusion of the Artinian algebras $B \hookrightarrow A$ where
\[
A=F[x_1 , x_2, \ldots, x_n]/(x_1^2, x_2^2, \ldots, x_n^2),
\]
\[
B=F[y_1, \ldots, y_r]/(y_1^{n_1+1}, y_2^{n_2+1}, \ldots, y_r^{n_r+1}).
\]
Recall that $l$ is the sum of the variables in $A$, so we may also write
$l=y_1 + y_2 + \cdots + y_r$. In particular, $l$ is an element of $B$.  There is a unique element of highest degree in $B$, namely
\[
l^n=n!x_1x_2\cdots x_n=\frac{n!}{n_1!\cdots n_r!} \ y_1^{n_1}\cdots y_r^{n_r}.
\]
\prskip

The inclusion $B \hookrightarrow A$ is actually a grade--preserving inclusion. Thus we have
the commutative diagram:
\prskip

\newcommand{\hookuparrow}{\mathrel{\rotatebox[origin=c]{90}{$\hookrightarrow$}}}
\newcommand{\hookdownarrow}{\mathrel{\rotatebox[origin=c]{-90}{$\hookrightarrow$}}}

\newcommand{\Ker}{{\rm Ker}}
\renewcommand{\Im}{{\rm Im}}
\newcommand {\str}{{ \stackrel{\times l}{\rightarrow} }}
\newcommand {\strtwo}{{ \stackrel{\times l^2}{\rightarrow} }}
\newcommand {\strc}{{ \stackrel{\times l^{c-2i}}{\rightarrow} }}

$
\begin{array}{cccccccccccccc}
A_0 & \str & A_1 & \str & A_2   & \str & \cdots &\str & A_{n-2} &\str & A_{n-1} & \str  & A_n \\
\hookuparrow & &\hookuparrow & &\hookuparrow &   & \cdots  &  & \hookuparrow &  &\hookuparrow  & &\hookuparrow  \\
B_0 &\str & B_1 &\str & B_2  & \str & \cdots  & \str & B_{n-2} & \str &B_{n-1}  &\str &  B_n
\end{array}
$
\prskip

\noindent In particular, for each $k=0,1, \cdots, [n/2]$, we have the commutative diagram
\[
\begin{array}{cccc}
\times l^{n-2k}:  &A_k & \to & A_{n-k} \\
   & \hookuparrow & & \hookuparrow \\
\times l^{n-2k}:& B_k  &  \to    &  B_{n-k}
\end{array}.
\]
Inasmuch as the map $\times l^{n-2k}:A_k \to A_{n-k}$ is  bijective, the restricted map $\times l^{n-2k}:B_k \to B_{n-k}$ is injective.  It is easy to see that $\dim _FB_k = \dim _F B_{n-k}$, whence the restricted map is also in fact bijective. Thus we have established the following theorem.
\prskip

\begin{theorem}
 Let $B$ be the monomial complete intersection ring,
\[
B=F[y_1, y_2, \cdots, y_r]/(y_1^{d_1+1}, y_2^{d_2 + 1}, \cdots, y_r^{d_r+1}).
\]
Then the multiplication map
\[
\times (y_1 + y_2 + \cdots + y_r)^{n-2k}: B_k \to B_{n-k}
 \]
is a bijection.
\end{theorem}
\pfskip
\prskip

\begin{remark}
{\rm(1)}
Ikeda has shown that if  a graded Artinian $F$-algebra $A$ has the strong Lefschetz property with a Lefschetz element   $l$, then $A[x] / (x^2)$  has the strong Lefschetz property with $x + l$ as a Lefschetz element. It immediately follows that the quadratic monomial complete intersection ring has the strong Lefschetz property with Lefschetz element $x + l$. Our results on the Weyr structures of blocked matrices can be used to prove this directly, in fact for $A[x]/(x^t)$.  This  enables us to determine the Weyr form for $\times l \in {\rm End} _F(B)$ for a general element $l$ in the monomial complete intersection ring~ $B$.\\
{\rm(2)} Let $A=F[x,y] / (x^m, y^n)$  be a monomial complete intersection ring in two variables. Assume that $F$ has characteristic  0 or $p > m + n - 2$. Our results on Weyr structures can be used to obtain the Weyr structure of the multiplication map $\times (x+y) : A \rightarrow A$, whereas this is not so easy to prove just by commutative algebra {\rm(}cf. \cite{II}{\rm)}. \hfill $\square$
\end{remark}
\prskip

\begin{theorem} \label{slp_monomial_ci}
Let $B$ be the monomial complete intersection ring,
\[
  B=F[x_1, \cdots, x_n]/(x_1^{d_1+1}, \cdots, x_n^{d_n+1}).
\]
Let $l$ be a general element in $B$.  Write $|B_i|$ for the dimension of $B_i$ (where $B_i$ is the homogeneous
space of $B$ of degree $i$). Then the Weyr structure of $\times l$ is given by the partition $|B|=|B_0|+|B_1|+ \cdots + |B_{d_1+ \cdots + d_n}|$, once arranged in decreasing order {\rm(}see also Proposition~\ref{multinomial_coef}{\rm)}.
\end{theorem}

\begin{proposition} \label{multinomial_coef}
Let $B$ be the monomial complete intersection ring as defined in Theorem~\ref{slp_monomial_ci}. Then the dimensions $\dim B_i$ of the homogeneous components of $B$ are determined as the coefficients of $T^i$ in the polynomial
\[
\prod _{j=1}^n(1+T+ T^2 + \cdots + T^{d_j}).
\]
\end{proposition}

\begin{remark}
Theorem~\ref{slp_monomial_ci} was proved by R.\ Stanley in \cite{STAN} using the Hard Lefschetz Theorem in algebraic geometry. It was also proved in \cite{WAT}, in which the second author used the theory of the Lie algebra $sl(2)$. It is rather an amazing fact that Theorem~\ref{slp_monomial_ci} is an easy consequence of Theorem~\ref{slp_quad_mono_ci},
since Theorem~\ref{slp_monomial_ci} is a generalisation of Theorem~\ref{slp_quad_mono_ci}. Theorem~\ref{slp_quad_mono_ci} was proved by Ikeda~\cite{IKE} in an elementary manner without a reference to the general case. Our proof here is elementary also, but quite different from the proof of Ikeda.
\end{remark}
\sectskip

\end{document}